\documentclass[11pt]{amsart}

\usepackage[T1]{fontenc}
\usepackage[utf8]{inputenc}
\usepackage{amsmath,amssymb,mathtools}
\usepackage{geometry}
\usepackage[colorlinks=true,linkcolor=red,citecolor=blue,urlcolor=blue]{hyperref}

\numberwithin{equation}{section}
\theoremstyle{plain}
\newtheorem{theorem}{Theorem}[section]
\newtheorem{lemma}[theorem]{Lemma}
\newtheorem{corollary}[theorem]{Corollary}
\newtheorem{prop}[theorem]{Proposition}
\theoremstyle{definition}

\theoremstyle{remark}

\newcommand{\M}{\mathbb{M}}

\newcommand{\U}[1]{\mathopen{\vert\!\vert\!\vert}#1\mathclose{\vert\!\vert\!\vert}}
\newcommand{\Kn}[2]{\lVert #1\rVert_{(#2)}}
\DeclareMathOperator{\Ree}{Re}

\title[Crossed-Heinz inequalities]{Crossed-Heinz matrix inequalities outside the unit interval}

\author{Jean-Christophe Bourin}
\address{Universit\'e Marie et Louis Pasteur, CNRS, LmB (UMR 6623),
F-25000 Besan\c{c}on, France}
\email{jean-christophe.bourin@umlp.fr}

\author{Mohammad Sababheh}
\address{Program of Mathematics, Abdullah Al Salem University,
Khaldiya, Kuwait City, Kuwait}
\address{Department of Basic Sciences, Princess Sumaya University for Technology,
Amman, Jordan}
\email{sababheh@yahoo.com}


\begin{document}

\begin{abstract}
In this paper, we consider the crossed-Heinz inequality
\[
\U{A^pB^q+B^pA^q}\le \U{A^{p+q}+B^{p+q}},
\]
where $A,B$ are positive definite matrices and $p,q>0$.
The validity of this inequality for every unitarily invariant norm
has been investigated since a question of the first author in 2009.
We prove the reverse inequality
\[
\U{A^{-p}B^q+B^{-p}A^q}
\ge \U{A^{q-p}+B^{q-p}},
\]
for every unitarily invariant norm, as a consequence of a stronger
eigenvalue-majorization inequality. Our main tool is the matrix
geometric mean.
\end{abstract}

\maketitle

	\section{Introduction and preliminaries}
	
	The singular values of $X\in\M_n$, the space of complex matrices of order $n$, are denoted by
	$s_1(X)\ge \cdots\ge s_n(X)$. 
	A norm $\|\cdot\|$ is unitarily invariant if
	$\U{UXV}=\U{X}$ for all $X\in\M_n$ and all unitary $U,V\in\M_n$.
	For $1\le k\le n$, the Ky Fan $k$-norm is 
	$\Kn{X}{k}=\sum_{j=1}^k s_j(X)$.
	For decreasing vectors $x,y\in\mathbb R^n$, we write
	$x\prec_w y$ when
	$\sum_{j=1}^k x_j\le\sum_{j=1}^k y_j$ for $1\le k\le n$.
	When $x$ and $y$ are vectors of singular values, this is equivalent to domination in every
	unitarily invariant norm. If $x,y$ are the vectors of eigenvalues, in decreasing order, of the Hermitian matrices $X,Y$, we also write $X\prec_w Y$, or $Y\succ_w X$, for $x\prec_w y$. This is the weak majorization order on the Hermitian matrices of $\M_n$. See the book \cite{BhatiaBook} for more on these classical notions.
	
	For positive definite $A,B\in\M_n$ and $t\in\mathbb R$, put
	\begin{equation*}
		H_t(A,B)=A^tB^{1-t}+A^{1-t}B^t,\qquad
		C_t(A,B)=A^tB^{1-t}+B^tA^{1-t}.
	\end{equation*}
	These are the Heinz and crossed-Heinz expressions, respectively. Note that the crossed-Heinz expression is symmetric in $A$ and $B$, a pleasant property for a mean.
	The classical Heinz inequality asserts that
	\begin{equation}\label{eq:heinz}
		\U{H_t(A,B)}\le\U{A+B},\qquad 0\le t\le1,
	\end{equation}
	for every unitarily invariant norm; see \cite{BhatiaBook,BhatiaDavis}.
	The function $t\mapsto\U{H_t(A,B)}$ is convex and symmetric
	about $t=1/2$; see \cite{BhatiaDavis}. As a consequence,
	\eqref{eq:heinz} is reversed when $t\notin[0,1]$.

	The question raised in \cite{Bourin2009} is whether
	\begin{equation}\label{eq:bourin}
		\U{C_t(A,B)}\le\U{A+B},\qquad 0<t<1,
	\end{equation}
	holds for all positive semidefinite matrices $A,B$ and every
	unitarily invariant norm; note that \eqref{eq:bourin} is equivalent to the first inequality in the abstract. Inequality \eqref{eq:bourin} is easy to prove
	 for the trace norm.

	Several partial results on this question
	were obtained in
	\cite{Bhatia2014, HayajnehEtAl2022b,HayajnehEtAl2023b}, and results in the setting of von Neumann algebras were given in \cite{LiuHeZhao,Zhang}.
	The comparison $\U{C_t(A,B)}\le\U{H_t(A,B)}$
	was shown to fail in general in \cite{Bottazzi}.
	The inequality \eqref{eq:bourin} was established in
	\cite{KittanehRicard} for $1/4\le t\le3/4$, in every dimension.
	In \cite{SababhehBourin}, we
	proved
	\begin{equation*}
		\U{C_t(A,B)}\le\U{H_t(A,B)}\le\U{A+B},\qquad
		A,B\in\M_2,\quad 0<t<1.
	\end{equation*}

	The operator norm case of \eqref{eq:bourin} was recently
	proved in all dimensions by Tian and Xu \cite{TianXu}, while the case of the Ky Fan $(n-1)$-norm on $\M_n$ was settled by Nicolas \cite{N}, thereby establishing \eqref{eq:bourin} in dimension three.
	The other intermediate Ky Fan norm cases in dimension $n\ge 4$
	remained open until very recently, when Tian and Xu \cite{TianXu2} announced a proof of \eqref{eq:bourin} in full generality.

	Unlike the ordinary Heinz expression, the crossed-Heinz
	expression does not necessarily have a convex norm as
	a function of the parameter, even for positive definite
	matrices of order two and the trace norm. Thus, the usual
	convexity argument for the Heinz inequality does not
	apply to the crossed-Heinz setting. Our approach instead
	relies on the matrix geometric mean and eigenvalue
	majorization.

	Motivated by the reverse of \eqref{eq:heinz}, we study the complementary
	question of whether
	\begin{equation}\label{eq:goal}
		\U{C_t(A,B)}\ge\U{A+B},\qquad t\notin[0,1].
	\end{equation}
Since \( C_t^*(A,B)=C_{1-t}(A,B)\), \eqref{eq:goal} is equivalent to
\begin{equation*}
		\U{C_{-t}(A,B)}\ge\U{A+B},\qquad t>0,
	\end{equation*}
which is clearly equivalent to the second inequality in the abstract when $q>p>0$. The case $p>q>0$ follows by taking adjoints and by replacing $A,B$ by their inverses. A limit argument yields the case $p=q>0$.

\section{Majorization outside the unit interval}

	For positive definite $A,B$ and $r\in\mathbb R$, set $P_r:=A^r+B^r$. Then
	\begin{equation*}
		C_r(A,B)=P_rP_{1-r}-(A+B).
	\end{equation*}
	This simple decomposition, involving three positive definite matrices, is the key to the proof below, in particular because $P_rP_{1-r}$ has the same eigenvalues as $P^{1/2}_rP_{1-r}P^{1/2}_r$. It also suggests taking the real part of $C_r(A,B)$ and looking for a domination of $A+B$ in the weak majorization order on Hermitian matrices. Another interesting feature of this decomposition is that $P_rP_{1-r}$ and 
$P^{1/2}_rP_{1-r}P^{1/2}_r$ look like geometric means, which suggests working with the matrix geometric mean $\#$, a powerful tool in matrix analysis since the seminal paper of Ando \cite{Ando}.

We briefly recall some basic facts about the matrix geometric mean and refer to the book \cite{BhatiaPDM} for an excellent account. Given two distinct positive definite matrices $A,B$, there is a path of positive definite matrices,
\begin{equation*}
t\mapsto A\#_t B=A^{1/2}(A^{-1/2}BA^{-1/2})^tA^{1/2}, \qquad t\in(-\infty,\infty).
\end{equation*}
Here $A\#_tB$ can be regarded as the point with barycentric coordinates $(1-t,t)$ with respect to the pair $\{A,B\}$ on this path (this path is in fact the constant-speed geodesic through $A$ and $B$ for a natural Riemannian manifold structure on the positive definite cone). In particular, $A\#_{1-t}B=B\#_tA$. The geometric mean $A\#_{1/2}B$ is simply denoted by $A\#B$; thus, for instance, $A\#(A\#B)=A\#_{1/4}B$. For $\theta\in[0,1]$, the homogeneous map
$$
(A,B)\mapsto A\#_{\theta}B
$$
is monotone and jointly concave (for the usual Loewner order $\le$ on Hermitian matrices). Thus, for two pairs of positive definite matrices
$(X_0,Y_0)$ and $(X_1,Y_1)$, we have, using homogeneity,
$$
(X_0+X_1)\#_{\theta} (Y_0+Y_1) \ge X_0\#_{\theta} Y_0 + X_1\#_{\theta} Y_1.
$$
This is a matrix H\"older inequality. Applying it to $X_0=A^x$, $Y_0=A^y$, $X_1=B^x$, $Y_1=B^y$, we see that our matrices $P_r=A^r+B^r$ satisfy a
H\"older-type inequality:

\begin{lemma}\label{holder} Let $A,B$ be positive definite. If $0\le \theta\le 1$, then for all real exponents $x,y$,
$$
(A^x+B^x)\#_{\theta}(A^y+B^y) \ge A^{(1-\theta)x +\theta y} + B^{(1-\theta)x +\theta y}.
$$
\end{lemma}

We also recall that the geometric mean of two positive definite matrices $A$ and $B$ is characterized as the unique positive definite matrix of the form $A^{1/2}UB^{1/2}$ with $U$ unitary. A classical log-majorization then ensures that
\begin{equation}\label{logmaj}
\log (A\#B) \prec_w \log |B^{1/2}A^{1/2}|,
\end{equation}
where $|Z|=(Z^*Z)^{1/2}$ is the right modulus of $Z\in\M_n$.

	\begin{theorem}\label{thm:main} Let $A,B\in\M_n$ be positive definite, let $t>0$ and set $C_{-t}:=C_{-t}(A,B)$. Then
\[
\frac{C_{-t}+C_{-t}^*}{2} \succ_w A+B.
\]		
	\end{theorem}

Note that \( C_{-t}+C_{-t}^*= C_{1+t}+C_{1+t}^*\).	

	\begin{proof} Since the functional $\Gamma_k(\cdot)$, the sum of the $k$ largest eigenvalues, is subadditive on the Hermitian matrices of $\M_n$, we have
\begin{equation*}
\Gamma_k\left(\frac{C_{-t}+C_{-t}^*}{2}\right) \ge \Gamma_k\left( \frac{P_{-t}P_{1+t}+P_{1+t}P_{-t}}{2}\right)-\Gamma_k(  A+B).
\end{equation*}
The Schur triangularization of \( P_{-t}P_{1+t} \) and the fact that the diagonal entries of a Hermitian matrix are majorized by its eigenvalues show that
 \begin{equation*}
\Gamma_k\left( \frac{P_{-t}P_{1+t}+P_{1+t}P_{-t}}{2}\right) \ge \Gamma_k\left(P_{-t}^{1/2}P_{1+t}P_{-t}^{1/2}\right) \ge \Gamma_k \left((P_{-t}\#P_{1+t})^2\right)
\end{equation*}	
where the second inequality is a standard consequence of the log-majorization \eqref{logmaj}.
Hence it suffices to show that
 \begin{equation} \label{step}
\Gamma_k \left((P_{-t}\#P_{1+t})^2\right) \ge 2\Gamma_k(A+B).
\end{equation}
Observing that $2(A+B)=(P_0\#P_1)^2$, we are led to compare $P_{-t}\#P_{1+t}$ and $P_0\#P_1$. Note that in the interval $[-t,1+t]$, the point $0$ has barycentric coordinates $(1-\theta,\theta)$ for some $\theta\in [0,1]$, while the point $1$ has barycentric coordinates $(\theta, 1-\theta)$. From Lemma \ref{holder} we infer
\[ 
P_0\le P_{-t}\#_{\theta} P_{1+t}, \qquad  P_1\le P_{-t}\#_{1-\theta} P_{1+t}.
\]
Therefore
\[
P_0\# P_1 \le \left( P_{-t}\#_{\theta} P_{1+t}\right)\#\left(  P_{-t}\#_{1-\theta} P_{1+t}\right) =  P_{-t}\# P_{1+t}.
\]
This shows that $(P_0\# P_1)^2\le V(P_{-t}\#P_{1+t})^2V^*$ for some unitary $V\in\M_n$, which is stronger than \eqref{step}.	\end{proof}

We derive several corollaries.
	
	\begin{corollary}\label{cor:crossed}
		 If $A,B\in\M_n$ are positive definite and $t<0$ or $t>1$, then
		\begin{equation*}
			A+B\prec_w\Ree C_{t}(A,B)
			\prec_w \left|\Ree C_{t}(A,B)\right|\prec_w | C_{t}(A,B)|.
		\end{equation*}
\end{corollary}

\begin{proof} The first weak majorization is Theorem \ref{thm:main}, the second one is trivial, and the third one is merely the norm inequality $\| X+X^*\| \le 2\|X\|$. \end{proof}

The weak majorization $A+B\prec_w | C_{t}(A,B)|$ is equivalent to the following norm inequality.

\begin{corollary}\label{cor:reverse}
		 If $A,B\in\M_n$ are positive definite and $t<0$ or $t>1$, then for every unitarily invariant norm,
				\begin{equation*}
			\U{A^tB^{1-t}+B^tA^{1-t}}\ge\U{A+B}.
		\end{equation*}
			\end{corollary}

As explained at the end of the introduction, this can be restated in the following form.

\begin{corollary}\label{cor:reverse2}
		 If $A,B\in\M_n$ are positive definite and $p,q>0$, then for every unitarily invariant norm,
				\begin{equation*}
			\U{A^{-p}B^{q}+B^{-p}A^{q}}\ge\U{A^{q-p}+B^{q-p}}.
		\end{equation*}
			\end{corollary}

In particular, taking $p=q=1$, we obtain the following, where $I$ denotes the identity matrix.

\begin{corollary}\label{cor:reverse3}
		 If $A,B\in\M_n$ are positive definite, then for every unitarily invariant norm,
				\begin{equation*}
			\U{A^{-1}B+B^{-1}A}\ge2\U{I}.
		\end{equation*}
			\end{corollary}

For decreasing vectors $x,y\in\mathbb{R}^n$ with positive entries, we write $x\prec_{w\log}y$ when
$\prod_{j=1}^k x_j\le\prod_{j=1}^k y_j$ for $1\le k\le n$. This weak log-majorization implies the weak majorization $x\prec_w y$; see \cite{BhatiaBook}. If $x,y$ are the vectors of singular values of $X,Y$, we also write $X\prec_{w\log} Y$ for $x\prec_{w\log}y$. Since the eigenvalues of $A^{-1}B$ lie in $(0,\infty)$, Corollary \ref{cor:reverse3} can be strengthened to a weak log-majorization, as a consequence of the following simple proposition.

\begin{prop}\label{prop:logmaj}
If $X\in\M_n$ has its spectrum lying in $(0,\infty)$, then
\[
I\prec_{w\log} \left| \frac{X+X^{-1}}{2}\right|.
\]
\end{prop}

\begin{proof} By Schur triangularization we may assume that
$X$ is upper triangular with diagonal part $D:=\mathrm{diag}(d_1,\ldots,d_n)$ where $d_j>0$ are the eigenvalues of $X$. Hence the diagonal part of the triangular matrix $(X+X^{-1})/2$  has all entries greater than or equal to $1$. A classical log-majorization of Weyl between the moduli of the eigenvalues and the singular values then completes the proof.
\end{proof}

\begin{corollary}\label{cor:reverse4}
 If $A,B\in\M_n$ are positive definite, then
\[			
I\prec_{w\log} \left| \frac{A^{-1}B+B^{-1}A}{2}\right|.
\]
\end{corollary}

 However, two singular values of $(A^{-1}B+B^{-1}A)/2$ may be smaller than $1$ when $n=3$. For instance, with
\[
A=\begin{bmatrix}1&3&6\\3&10&22\\6&22&53\end{bmatrix},\qquad
B=\begin{bmatrix}4&12&24\\12&38&80\\24&80&177\end{bmatrix},
\]
we have
\[
 \frac{A^{-1}B+B^{-1}A}{2}=\begin{bmatrix}17/8&21/8&15/4\\0&5/4&1\\0&0&1\end{bmatrix},
\]
whose singular values are approximately $5.293$, $0.773$ and $0.650$.

\section{The operator norm in \texorpdfstring{$C^*$}{C*}-algebras}

For the operator norm, that is for $k=1$, the proof of Theorem \ref{thm:main} works in any unital $C^*$-algebra $\mathcal{A}$. Positive invertible elements of $\mathcal{A}$ are again called positive definite, and $C_t(A,B)$, $P_r$ and $A\#_\theta B$ are defined by the continuous functional calculus as in Section 2. Representing $\mathcal{A}$ faithfully on a Hilbert space, we see that the monotonicity and the joint concavity of $\#_\theta$ \cite{KuboAndo}, hence Lemma \ref{holder}, remain valid in $\mathcal{A}$, as well as the identity $(X\#_\theta Y)\#(X\#_{1-\theta} Y)=X\#Y$. We denote by $\|\cdot\|_\infty$ the norm of $\mathcal{A}$, and for a Hermitian $h\in\mathcal{A}$, we set $\lambda_{\max}(h)=\max\sigma(h)=\sup_\varphi\varphi(h)$, where $\varphi$ runs over the states of $\mathcal{A}$.

\begin{theorem}\label{thm:cstar} Let $\mathcal{A}$ be a unital $C^*$-algebra, let $A,B\in\mathcal{A}$ be positive definite, and let $t<0$ or $t>1$. Then
\[
\lambda_{\max}\left(\Ree C_t(A,B)\right)\ge\|A+B\|_\infty.
\]
In particular, $\|A^tB^{1-t}+B^tA^{1-t}\|_\infty\ge\|A+B\|_\infty$.
\end{theorem}

\begin{proof} Since $C_t^*(A,B)=C_{1-t}(A,B)$, it suffices to show that $\lambda_{\max}(\Ree C_{-t})\ge\|A+B\|_\infty$ for $t>0$, where $C_{-t}:=C_{-t}(A,B)$. Set $P=P_{-t}$ and $Q=P_{1+t}$. Since $\lambda_{\max}$ is subadditive,
\[
\lambda_{\max}(\Ree C_{-t})\ge\lambda_{\max}(\Ree PQ)-\|A+B\|_\infty.
\]
The spectrum of an element of $\mathcal{A}$ is contained in its
numerical range $\{\varphi(\cdot)\}_\varphi$. Moreover,
\[
PQ=P^{1/2}(P^{1/2}QP^{1/2})P^{-1/2}.
\]
Hence
\[
\lambda_{\max}(\Ree PQ)=\sup_\varphi\Ree\varphi(PQ)\ge\max\sigma\left(P^{1/2}QP^{1/2}\right)=\left\|P^{1/2}QP^{1/2}\right\|_\infty=:\mu.
\]
Now $P^{1/2}QP^{1/2}\le\mu I$, that is $Q\le\mu P^{-1}$, so that $P\#Q\le P\#(\mu P^{-1})=\mu^{1/2}I$. Thus $\mu\ge\|P\#Q\|_\infty^2$. Finally, the barycentric argument of the proof of Theorem \ref{thm:main} yields $P\#Q\ge P_0\#P_1=(2(A+B))^{1/2}$, hence $\|P\#Q\|_\infty^2\ge2\|A+B\|_\infty$. Combining these inequalities gives $\lambda_{\max}(\Ree C_{-t})\ge\|A+B\|_\infty$. The last assertion follows from $\|X\|_\infty\ge\|\Ree X\|_\infty\ge\lambda_{\max}(\Ree X)$.
\end{proof}

As in the introduction, the norm inequality in Theorem \ref{thm:cstar} can be restated as $\|A^{-p}B^q+B^{-p}A^q\|_\infty\ge\|A^{q-p}+B^{q-p}\|_\infty$ for all $p,q>0$. Applying it to the Calkin algebra, we obtain:

\begin{corollary}\label{cor:ess} Let $A,B$ be positive invertible operators on an infinite-dimensional Hilbert space. If $t<0$ or $t>1$, then
\[
\|A^tB^{1-t}+B^tA^{1-t}\|_{\rm ess}\ge\|A+B\|_{\rm ess}.
\]
\end{corollary}

\begin{proof} Let $\pi$ be the quotient map onto the Calkin algebra, so that $\|X\|_{\rm ess}=\|\pi(X)\|_\infty$. Since $\pi$ is a unital $*$-homomorphism, $\pi(A)$ and $\pi(B)$ are positive definite and $\pi(A^r)=\pi(A)^r$, $\pi(B^r)=\pi(B)^r$ for all $r\in\mathbb R$. Hence $\pi(C_t(A,B))=C_t(\pi(A),\pi(B))$, and Theorem \ref{thm:cstar} completes the proof.
\end{proof}

	\section*{Declarations}
\noindent\textbf{Conflict of interest.}
The authors declare that there is no conflict of interest.
\medskip

\noindent\textbf{Funding.}
The authors received no specific funding for this work.
\medskip

\noindent\textbf{Data availability.}
No datasets were generated or analyzed during the current study.
\medskip

\noindent\textbf{Author contributions.}
The authors contributed to the conception of the work, the development
of the mathematical results, the proofs, and the preparation of the manuscript.

\section*{Declaration of generative AI and AI-assisted technologies
in the writing process}
The authors used OpenAI ChatGPT and Anthropic Claude for language editing,
\LaTeX{} formatting, and preliminary consistency checks. The authors retained
full control over the mathematical development and presentation of the results
and take full responsibility for the content of the publication.

% The bibliography is alphabetized by the first author's surname.
	
\end{document}